\documentclass[hidelinks,onefignum,onetabnum]{siamart220329}

\usepackage{amsfonts}
\usepackage{graphicx}
\headers{Sigmoid functions and multiscale resolution}{Daan Huybrechs and Lloyd N. Trefethen}

\title{Sigmoid functions and multiscale resolution of singularities
\thanks{Submitted to the editors DATE.}}

\author{Daan Huybrechs\thanks{Dept.\ of Computer Science, KU Leuven, 3001 Leuven, Belgium}
  (\email{daan.huybrechs@kuleuven.be}).
\and Lloyd N. Trefethen\thanks{Mathematical Institute, University of Oxford, Oxford OX4 4DY, UK}
  (\email{trefethen@maths.ox.ac.uk}).}

\def\Im{\hbox{Im\kern .5pt}}
\def\intz{[\kern .4pt 0,1]}
\def\zint{[-\infty,0\kern .4pt]}
\def\pdes{PDE\kern .4pt s}
\def\e{\varepsilon}
\def\ak{a_k^{}}\def\ek{\e_k^{}}\def\sk{s_k^{}}\def\sj{s_j^{}}
\def\rn{r_n^{}}\def\Rn{R_n^{}}
\def\smin{s_{\min{}}}

\begin{document}

\maketitle

\begin{abstract}
In this short, conceptual paper we observe that essentially the same
mathematics applies in three contexts with disparate literatures:
(1) sigmoidal and RBF approximation of smooth functions,
(2) rational approximation of analytic functions near singularities, and
(3) $hp$ mesh refinement for solution of \pdes.
The relationship of (1) and (2) is as simple as the change of variables $s = \log(x)$,
and our informal mnemonic for this relationship is ``sigmoid = log(ratapprox).''
\end{abstract}



\section{\boldmath Sigmoids $\leftrightarrow$ rational approximation}
\begin{figure}
\vskip 5pt
\begin{center}
\includegraphics[scale=.95]{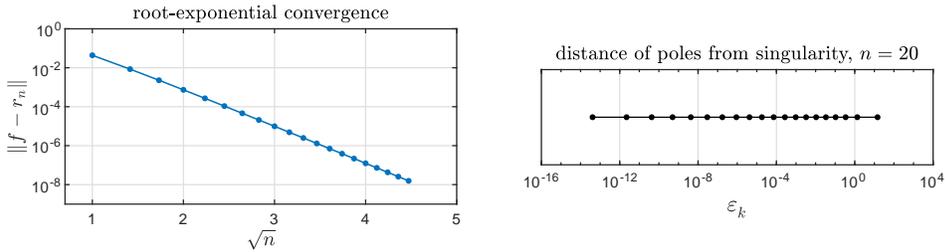}~~
\end{center}
\vskip -5pt
\caption{\label{fig1}Root-exponential convergence (\kern 1pt left, $0\le n \le 20$) and
exponential clustering of poles near\/ $0$ in $(-\infty,0)$
(right, $n=20$) for degree $n$ minimax rational approximation
of\/ $f(x) = \sqrt x$ on $\intz$.  Note that the horizontal
axis in the first plot is $\sqrt n$.  These effects generalize to rational
approximation near any branch point singularity.}
\end{figure}

Functions with branch point singularities can be approximated with root-exponential
convergence by rational functions.  For example, consider
\begin{equation}
f(x) = \sqrt x, \quad x\in \intz.
\label{modprob}
\end{equation}
Since Donald Newman in 1964~\cite{newman} it has been known that there are degree
$n$ rational functions
\begin{equation}
\rn(x) = a_0^{} + \sum_{k=1}^n {\ak \over 1 + x/\ek} 
\label{partfrac}
\end{equation}
such that 
\begin{equation}
\|f - \rn\kern .7pt  \| =  O(\exp(-C \sqrt n\kern 1.5pt)), \quad C>0,
\label{conv}
\end{equation}
where $\|\cdot\|$ is the supremum norm on $\intz$.  What makes
this {\em root-exponential convergence\/} possible is that
the poles $\{-\ek\}$ are exponentially clustered near the
singularity at $x=0$.  This effect applies near any branch point
singularity~\cite{lightning,herremans,clustering} of a real or
complex function and is the basis of ``lightning PDE solvers''
for the Laplace, biharmonic, and Helmholtz equations in domains
with corners~\cite{baddoo,stokes,lightning,helm}.  Figure~\ref{fig1}
illustrates root-exponential convergence and exponential clustering
for this model problem.

The function $1/(1+x/\ek)$ of (\ref{partfrac}) is monotonically decreasing
for $x\in\intz$ and takes values ${\approx}\, 1$ for $x\ll \ek$ and
${\approx}\,0$ for $x\gg \ek$.
If we introduce the change of variables
\begin{equation}
s = \log x \in \zint, \quad x = e^s \in \intz,
\label{xs}
\end{equation}
then this function is transformed into
\begin{equation}
{1\over 1+e^s/\ek} =
{1\over 1+e^{s-\sk}},
\label{act}
\end{equation}
where $\sk = \log \ek$.  This function, or more properly its
reverse $1/(1+e^{\sk-s}),$ is the most basic example of a sigmoid
function, known as the {\em logistic function}.  In physics
it goes by the name of the {\em Fermi\/} or {\em Fermi--Dirac
function,} and it is an elementary transformation of the hyperbolic
tangent.   Functions of this kind are prototypical activation
functions in neural networks, and the literature of this area is
vast~\cite{cybenko,goodfellow,hornik,lecun,pinkus}.

\begin{figure}
\vskip 5pt
\begin{center}
\includegraphics[scale=.95]{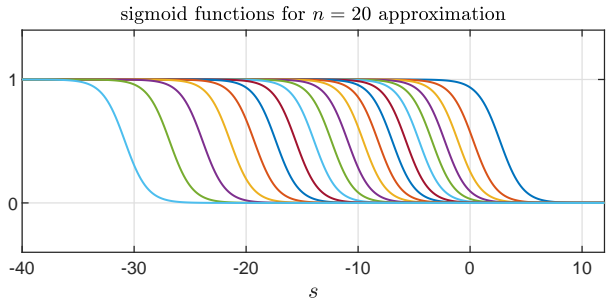}
\end{center}
\vskip -5pt
\caption{\label{fig2}In the $s = \log x $ variable, the exponentially clustered
poles of Figure~$\ref{fig1}$ become sigmoid functions $(\ref{act})$ translated to 
various center points $\sk$.  Root-exponential
convergence of rational approximations becomes a statement about approximation
of smooth functions by linear combinations of
translates of a fixed smooth function.}
\end{figure}

With the change of variables (\ref{xs}), we may follow
(\ref{modprob}) and (\ref{partfrac}) to define
\begin{equation}
F(s) = f(x) = e^{s/2}
\label{Ff}
\end{equation}
and
\begin{equation}
\Rn(s) = \rn(x) = a_0^{} + \sum_{k=1}^n {\ak \over 1 + e^{s-\sk}} .
\label{Rr}
\end{equation}
Equation (\ref{conv}) then implies that there are approximations
(\ref{Rr}) such that
\begin{equation}
\|F - \Rn \| =  O(\exp(-C \sqrt n\kern 1.5pt)), \quad C>0,
\label{conv2}
\end{equation}
where $\|\cdot\|$ is now the supremum norm on $\zint$.  Equation
(\ref{conv2}) is nothing else than a claim about root-exponential
approximation of a smooth function on $\zint$ by linear combinations
of translates of a standard smooth function, plus a constant.
Following results such as those of~\cite{driscoll}, this could
probably be proved more simply in this setting than in the setting
of rational functions, where the standard proof technique is the
relatively advanced Hermite integral formula~\cite{analcont}.
The root-exponential rate results from balancing discretization
errors associated with separations $\Delta \sk = O(1/\sqrt n \kern
1.4pt)$ against truncation errors associated with a grid extent
$\smin = \min \sk = O(\sqrt n \kern 1.4pt)$~\cite{herremans}.
Figure~\ref{fig2} illustrates the smooth functions in question for
the problem of Figure~\ref{fig1}.

To summarize this opening section: with the change of variables
$s = \log x$, the approximation of a smooth function by linear
combinations of translates of a logistic function becomes equivalent
to the approximation of a function with a branch point singularity
by rational functions with exponentially clustered poles.

\section{\boldmath Rational approximation $\leftrightarrow$ $hp$ mesh refinement}
Rational approximation with poles exponentially clustered
near singularities seems obviously akin to the resolution of
functions near singularities by exponentially refined meshes.
Such techniques of {\em mesh refinement\/} are well known in the
literature of the finite element method and associated approximation
theory~\cite{devore,guibab,guobab,melenk,scherer,schwab}.

Exponential clustering of poles is reflected in the approximately
uniform spacing on the semilogx scale in the right image of
Figure~\ref{fig1}, or equivalently, the approximately uniform
spacing of the sigmoid functions (\ref{act}) in Figure~\ref{fig2}.
However, it is notable that in both of these images the spacing
is only approximately uniform, growing sparser toward the left.
This is the phenomenon of {\em tapered exponential clustering\/}
investigated in~\cite{clustering}.  Quantitatively, one finds that
the density of poles with respect to the $s$ variable decreases about
linearly as $s$ decreases to some value $\smin$.  This distribution
brings a factor of $2$ improvement in convergence rate as a function
of $n$---because a uniform density would have the same convergence
rate but twice as many poles.  (The more local sparsification
in the rightmost few points of Figures~\ref{fig1} and~\ref{fig2}
is investigated in~\cite{herremans} with appeal to the asymptotic
results of Stahl~\cite{stahl}.)

We have noticed that this tapered exponential clustering corresponds
closely to what is known as $hp$ mesh refinement ($h$ stands for grid
spacing, $p$ for order of approximation).  In particular, the
standard $hp$ mesh refinement formula in one dimension has the
same linear pattern just described, with polynomial order taking
the role of pole density.  A singular function such as $\sqrt x$
on $\intz$ is approximated by piecewise polynomials on intervals
of lengths decreasing exponentially toward the singularity, with
polynomial representations of linearly decreasing degrees $\dots,
3,2,1,0.$ It is the same pattern, and it brings the same factor of
2 speedup for the same reason.

One can explain linear tapering in various ways in various settings.
In~\cite{clustering} an argument is given based on potential
theory.  Here is an outline of the simpler argument that originates with DeVore and
Scherer~\cite{devore,scherer} in the study of piecewise polynomial
approximations of $x^\alpha$ on exponentially graded meshes on $\intz$.
We speak for simplicity in terms of mesh refinement by factors of $1/2$,
though the optimal factor is actually $(\sqrt 2 - 1)^2 \approx 0.172$.
The observations of DeVore and Scherer were generalized
to ODE and PDE discretizations a few years later by Babu\v ska and his
collaborators~\cite{guibab,guobab}.  In multiple dimensions, the details
change.
\medskip

\begin{enumerate}

\item Approximation of $\sqrt x$ on $[1/2,1]$ is the same as
approximation of $\sqrt{2x}$ on $[1/4,1/2]$.

\item Therefore approximation of $\sqrt x$ on $[1/4,1/2]$ is the
same problem too, but with an accuracy criterion loosened by a
factor $\sqrt 2$.

\item Functions like these, bounded away from singularities, can
be approximated by polynomials with exponential convergence.

\item Therefore, that loosening by the factor $\sqrt 2$ allows one
to lower the degree of the polynomial by a constant increment and
still get the same accuracy.

\item Repeat on $[1/8,1/4]$, $[1/16,1/8],\dots.$

\end{enumerate}
\medskip 

To summarize this second section: the standard formula for $hp$ mesh
refinement in 1D involves a linear decrease of polynomial degree
toward the singularity, and this corresponds to the linear decrease
of pole density on a logarithmic scale in tapered exponential
clustering of poles in rational approximation, resulting in the
same factor of 2 speedup.

We note that this standard $hp$ mesh refinement strategy is not
the only way to achieve linear tapering.  An alternative would
be to hold the polynomial degree $p$ fixed and instead refine $h$
super-exponentially at the singularity.  We do not know if such a
prescription has been used in finite element calculations.

\section{Double exponential and generalized Gauss quadrature}
In the area of quadrature or numerical integration, many
methods have been developed for dealing with singularities.
When an endpoint singularity like $x^\alpha$ is known, a
targeted quadrature formula can be derived: the prototype is
Gauss--Jacobi quadrature.  For dealing with more complicated or
unknown singularities, however, more general techniques have been
proposed.  One is {\em double exponential\/} or {\em tanh-sinh\/}
quadrature~\cite{baileyb,ms01,tm74}.  As illustrated in Figure~14
of~\cite{clustering}, the tanh-sinh formula with standard parameter
choices produces a tapered exponentially clustered distribution of
quadrature points, indicating that it is probably related to what
is seen with rational approximation and $hp$ mesh refinement.

Kirill Serkh (private communication) has shown us that 
similar effects also arise with {\em generalized Gauss\/} and
{\em universal\/} quadrature formulas~\cite{bgr,brs}.  These are
quadrature formulas that are constructed by linear algebra methods
related to Gauss quadrature so as to be efficient at integrating not
just a single singularity such as a fixed power $x^\alpha$ but
a range of singularities such as $x^\alpha$, $\alpha \in \intz$.
Again it appears that in important cases, the nodes are exponentially
clustered near the singularity with a tapered distribution.

\section{\boldmath The physics of $s = \log(x)$: separation of scales}
The change of variables $s = \log x$ is not just an algebraic trick.
It also has a physical interpretation alluded to in section~5
of~\cite{clustering}; see for example Figure~11 of that paper.

Even before introducing the change of variables, the point
can be seen in the $x$ variable.  The function $1/(1+x/\ek)$, with its pole
at distance $\ek$ to the left of $x=0$, is essentially constant and hence inactive
to the right of $x=0$ for $x\ll \ek$ (taking the value $1$) and $x\gg \ek$ (taking the value
$0$).  It is only for $x\approx \ek$ that this function is active.  Thus
exponentially separated
poles $\{-\ek\} \subseteq (-\infty, 0)$ are physically decoupled, operating in
independent regimes, with each pole at $-\ek < 0$ affecting the approximation
on $(0,\infty)$ only for $x\approx \ek$.

\begin{figure}
\vskip 13pt
\begin{center}
\includegraphics[scale=.95]{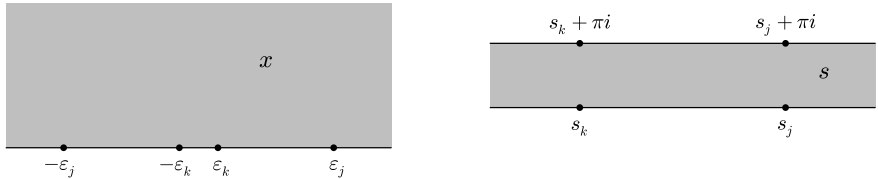}
\end{center}
\vskip -2pt
\caption{\label{fig3}Sketch of the ``physics'' of the change of
variables $s=\log x$.  This is a conformal map of the upper half $x$-plane to
an infinite strip in the $s$-plane, with poles $-\ek$ and
sample locations $\ek$ on exponentially separated scales mapping to
poles $\sk+\pi i$ and sample locations $\sk$ on opposite sides of the strip.
The well-known exponential decay of influences along strips explains why
a pole at $\sk+\pi i$ has a significant effect at $\sk$ but a much smaller
effect at a different point $\sj$ with $|\sj-\sk|\gg 1$.
}
\end{figure}

The change of variables $s = \log x$ suggests a physical explanation
of this separation of scales effect.  A problem with a singularity
at $x=0$ can be motivated as a model of a corner singularity in
a PDE problem.  Specifically, suppose a Laplace problem is posed
in the upper half complex $x$-plane with a singularity at $x=0$.
Changing to $s = \log x$ transplants this problem to the infinite
strip $0 < \Im s < \pi$ in the $s$-plane.  The problem is now smooth,
with the singularity moved to $-\infty$.  Exponentially clustered
poles $-\ek\in (-\infty,0)$ become well separated poles $\sk +
\pi i$ on the upper side of the strip.  And now, as sketched in
Figure~\ref{fig3}, it is a well-known effect of potential theory
(or elasticity, where it is called the St.\kern 1pt Venant principle)
that influences decay exponentially with distance along a strip.

The argument just made is tied to the Laplace equation, because its
solutions are invariant with respect to conformal maps.  However, the
essence of the matter will be the same for any problem whose highest
order derivative is the Laplacian, because close to a singularity,
this term will dominate.  With the Helmholtz equation $\Delta u +
k^2 u = 0$, for example, the influence of the $k^2 u$ term quickly
shuts off to zero as one comes exponentially close to a corner.

\section{Radial basis functions and other activation functions}
A rational function $r(x)$ as in (\ref{partfrac}) is a sum of simple
poles, which in the $s$ variable becomes a linear combination
of sigmoids as in Figure~\ref{fig2}.  The picture changes little
for various other activation functions.  For example, if the poles
$(x+\ek)^{-1}$ in (\ref{partfrac}) are replaced by powers
$(x+\ek)^{-a}$ for an arbitrary $a>0$, Figure~\ref{fig2} does not change very much.
This matches theoretical and experimental results in neural networks,
where choices between activation functions are typically based more on
the efficiency of learning algorithms such as stochastic gradient descent
than on approximation power.

As many authors have noted, a closely related topic is that of
approximation by radial basis functions.
Here again one approximates a complicated function by a linear combination of
translates of a simpler fixed function, and convergence may be very fast when
the latter is smooth~\cite{buhmann,driscoll,fasshauer,mhaskar,mmicc,park}.

In closing we note that although sigmoidal and other activation 
functions are important in neural networks and deep learning,
the present paper touches only the surface of that discipline. 
The approximation (\ref{Rr}) is not composite, but involves just what is conventionally
called a single ``hidden layer.''  Smooth activation functions have
accuracy advantages for single-layer approximation, but in the multi-layer
setting of deep learning, that advantage diminishes and
the simpler non-smooth function known as ReLU is
used more often~\cite{lecun,strang}.  We cannot resist mentioning that the universal
approximation power of ReLU units was exploited by Henri Lebesgue at age 23 in his first
published paper, which presented a new proof of the Weierstrass
approximation theorem~\cite{lebesgue}.

\section*{Acknowledgments}
We are grateful to Kirill Serkh of the University of Toronto for
showing us the tapered exponentially clustered nodes of generalized
Gauss and universal quadrature formulas.

\indent~~
\end{document}